\documentstyle[12pt]{article}

\catcode`@=11
\long
\def\@makefntext#1{\noindent #1}
\newskip\tabcentering \tabcentering=1000pt plus 1000pt minus 1000pt

\def\MCH#1#2{\setbox0=\hbox{\raise#1\hbox{#2}}\smash{\box0}}

\newcommand{\beqs}{\begin{eqnarray*}}
\newcommand{\eeqs}{\end{eqnarray*}}
\newcommand{\beqn}{\begin{eqnarray}}
\newcommand{\eeqn} {\end{eqnarray}}

\def\beq{\begin{equation}}
\def\eeq{\end{equation}}

\begin{document}
\title{\Large Long-time existence of mean curvature flow with
external force fields
\thanks{
 Supported by Natural Science Foundation of China (No. 10631020) and Basic Research  Grant of
 Tsinghua University
 }}

\author{Yannan  Liu and Huaiyu Jian \\
 {   Department of Mathematical Sciences}\\
 { Tsinghua
University,  Beijing 100084, P. R. China}\\
{\small Correspondence should be addressed to  Huaiyu Jian
(email:hjian@math.tsinghua.edu.cn)}\\ }
\date{}

\maketitle

{\bf Abstract} In this paper, we study the evolution of submannifold
 moving by mean curvature minus a external force field. We prove
 that the flow has a long-time smooth solution for all time  under
 almost optimal conditions.
  Those conditions are that
 the second fundamental form on the initial submanifolds is not too
 large, the external force field, with its any order
derivatives, is bounded, and    the field is convex with its
eigenvalues satisfying a pinch  inequality.  \\

 {\bf Key words } parabolic equation,
 mean curvature flow,  maximum principle (for tensor).

\newpage
\section{Introduction}
{\subsection {Background}}\indent In this paper, we study the flow
\begin{equation}\frac{dF}{dt} = - (H_{\alpha} -\omega_{\alpha}
){e_\alpha}\equiv -f_{\alpha}e_{\alpha}\label{eq1.1}\end{equation}
where
$$F_t:=F(\cdot,t): M^{n}\rightarrow R^{n+k}$$ is a family of smooth
immersions with  $M_t= F_{t}(M)$ and M is compact oriented
submanifold in $R^{n+k}$, H denotes the mean curvature vector of
$M_t$ w.r.t  unit normal field $e_{\alpha}$, $\alpha=n+1,\ldots,n+k$
, $\omega $ is a given smooth function in $R^{n+k}$,
$\overline{\nabla}\omega$ is the standard gradient field of
$\omega$ in
 $R^{n+k}$ , and  $\omega_{\alpha}=\langle\overline{\nabla}\omega,e_{\alpha}\rangle$
 .\\
\indent  This flow generalizes the well-known {\sl mean curvature
flow}, i.e., the case of $\omega \equiv  const $, and it comes
directly from the study of the Ginzburg-Landau vortex. As was
 shown in [1,2], there are two models which are, respectively,
 reduced to the Ginzburg-Landau system of parabolic equations
 \begin{equation}\frac{\partial V_{\varepsilon}}{\partial t}   =
\overline{\Delta }V_{\varepsilon} +\overline{\nabla} \omega
\overline{\nabla} V_{\varepsilon} +AV_{\varepsilon}
  + \frac{BV_{\varepsilon}}{\varepsilon^2}
(1-|V_{\varepsilon}|^2)\end{equation}in $R^{m}\times{(0,\infty)}$,
where $\varepsilon$ is a small positive parameter and $\omega,$ A,B,
are known functions. One  is a simple equation simulating
inhomogeneous type $\Pi$ superconducting materials [3], and the
other is a three-dimensional superconducting thin films having
variable thickness [4]. An important problem in Ginzburg-Landau
superconductors is to study the vortex dynamics, i.e, the
convergence of $V_\varepsilon$ as well as of their zero points
(which, roughly, are called vortex) as
$\varepsilon\rightarrow0$.\\
\indent When $m=2$ and the initial vortex consists of finite
isolated points, it was proved that the vortex dynamics of the
Dirichlet problem for (1.2) is described by the ODE system[1,5,6]:
$$ \frac{\partial x}{\partial t}= - \overline{\nabla} \omega(x).$$
When $m\geq 2$ and the initial vortex consists of a filment or even
a codimension k submanifold, it was proved [2] that as
$\varepsilon\rightarrow0$, the vortex of Cauchy problem for (1.2) is
approximated by the evolution of the initial vortex according to
flow (1.1) on the time internal in which the flow is smooth. Similar
results were obtained for Neumann problem in [7] and for case of
$\overline{\nabla} \omega= 0$ in [7,8].\\
\indent Therefore, it is  important in physics to consider the
long-time existence of the flow (1.1).

On the other hand, {\sl mean curvature flow} has been strongly
studied in last decades. It is well-known that the flow must blow up
in finite time except that the initial submanifolds are graphic, see
[9-16] for the details. Hence, it is natural to ask for what kind of
functions $\omega $ (1.1) has long-time existence.

 {\subsection {Main results}}\indent Higher co-dimension {\sl mean
 curvature flow}, i.e., (1.1) without external force field, has been studied
 in[9-13]], while there are a lot of studies on {\sl mean
 curvature flow}   for   hypersurfaces, see  [14][15][16] for
 example. All those papers show that {\sl mean
 curvature flow} must blow up and so singularity happens
in finite time, except that the initial surfaces  are entire graphs
or graphic submanifolds.
 In this
 paper, we are concentrated on the  long-time existence of (1.1). Here is the main
 results.

\indent {\bf Theorem 1.1}\indent {\sl If there exist positive
constants $C, C_3, \overline{\lambda} , \underline{\lambda} $ with
$\overline{\lambda} < 2\underline{\lambda} $ such that
 the following conditions are satisfied: \\
 (1) $\underline{\lambda}|\xi|^{2}\leq
\overline{\nabla}^{2}{\omega(x)}\xi_{i}\xi_{j}\leq
\overline{\lambda}|\xi|^{2} $    and
$|\overline{\nabla}^3\omega(x)|\leq  C_3$  for all $\xi \in R^{n+1}$
and for $x\in M_t,$ where $M_t$ is any solution of (1.1) on any
finite  time interval $[0, T];
$\\
 (2)  $|A|^{2} < C$ on $M_{0};$\\
 (3) there exist a $\delta >0$ such that $5a^{4} +
 |a|C_3 + (\overline{\lambda}-2\underline{\lambda})a^{2}\leq 0$  for
 all $a$ satisfying $\sqrt{C}-\delta <a<\sqrt{C};$\\
(4) $|\overline{\nabla}^{i}\omega(x)|$ is uniformly bounded for all
$x\in M_t$ and for $i=1, 2, 3, \cdots,$  where $M_t$ is any solution
of (1.1) on any
finite  time interval $[0, T];$ \\
 then the flow (1.1) has a smooth solution   for all time $t\in [0, \infty
 ).$\\ }
Throughout this paper, flow (1.1) is denoted by $(1.1)^{'}$ in the
 case of $k=1,$ i.e., the hypersurfaces case.

\indent {\bf Theorem 1.2}\indent {\sl Suppose that the assumptions
of theorem 1.1 are satisfies except that (1.1) is replaced by
$(1.1)^{'}$ and (3) is replaced by\\
$(3)^{'}$  there exist a $\delta >0$ such that $a^{4} +
 |a|C_3 + (\overline{\lambda}-2\underline{\lambda})a^{2}\leq 0$ for
 all $a$ satisfying $\sqrt{C}-\delta <a<\sqrt{C}+\delta;$ \\
 then the flow $(1.1)^{'}$ has a smooth solution   for all time $t\in [0, \infty
 ).$\\ }
\indent {\bf Remark 1.3}\indent An easier verification of
assumptions (1) and (4) is to assume that  they hold for all $x\in
R^{n+k}.$

The following theorem generalizes the convexity-preserved of mean
curvature flow in [14].

 \indent {\bf Theorem 1.4}\indent {\sl Let
$T>0$ and $M_t$ be a smooth solution of flow $(1.1)^{'}$ on the time
interval $[0, T].$ If $ \overline{\nabla}^{3}\omega \equiv 0$ and
$M_0$ is convex, then $M_t$ is convex for all $t\in [0, T].$}

 \indent  Physically,  $\omega$  is a  density function  and actually
has the form
 $$\omega=\frac{1}{2}(c_{1}x_{1}^{2}+ \ldots +
 c_{n+1}x_{n+1}^2)$$ for $c_i>0, $ ( see [3] for
 example),  but theorems 1.1 and 1.2 can not be applied directly to this special case,
 because
$|\overline{\nabla}\omega|$ is not known to be  bounded uniformly at
this moment. However, we can give the long-time existence for this
$\omega$ under hypersurfaces
case.\\
\indent {\bf Corollary 1.5}\indent {\sl Suppose that
$\omega=\frac{1}{2}(c_{1}x_{1}^{2}+ \ldots +
 c_{n+1}x_{n+1}^2)$  where $c_i$ are positive constants and let
 $M=\max c_i $ and $m=\min c_i .$ If $M<2m $ and $|A|^2 <2m-M$ on
 $M_0,$ then for any $T>0$ the flow $(1.1)^{'}$ has a smooth
 solution for all $t\in [0, T].$ }

 We would like to point out that Corollary 1.5 generalizes
 theorem 1.3 in [17] which studies  $(1.1)^{'}$ for the special case
 $\omega =c |x|^2.$ That theorem also shows that flow $(1.1)^{'}$
 must blow up in finite time  either if $c<0,$ or if $c>0$ and $|A|^2>c$ on $M_0,$
 which means that the both
 convexity of $\omega $ (as in assumption (1)) and the small of the initial $|A|^2$
 (as in assumptions (2) and (3) ) are necessary.

 \indent In section 2 we will give  notations and preliminaries , and we will
 give the  proofs of theorems 1.2, 1.4 and corollary 1.5 in section 3
 and the proof of theorem 1.1 in section 4.\\

\section{Preliminaries}
\setcounter{equation}{00} \indent Throughout this paper, we use the
following notations. $\langle,\rangle$ denotes the usual inner
product in $R^{n+k}$. If $M$ is given as in section 1 and  $F$
denotes its parametrization in $R^{n+k},$   the metric $\{g_{ij}\}$
are given by
$$ g_{ij}(x)=\langle\frac{\partial
F(x)}{\partial{x_i}},\frac{\partial F(x)}{\partial{x_j}}\rangle,
\indent x\in M.$$ Let $\nabla$ denotes the Levi-Civita connection on
M, while $\overline{\nabla}$ denotes the standard gradient in
    $R^{n+k}.$ We will use $i,j,k,\cdots $ to
 denote the tangent indexes  and $\alpha,\beta,\gamma,\cdots$ for normal ones.
 Doubled indices  always mean to sum from 1 to $n$ for $i,j,k,\cdots $
 and from 1 to $k$ for $\alpha,\beta,\gamma,\cdots .$
  Indices are raised and lowered w.r.t $g^{ij}$ and
 $g_{ij}.$
Moreover, we will identify $V\in{T_{x}M}$ with
    $DF(V)\in{T_{F(x)}R^{n+k}}$. Also, we will use
    $\langle,\rangle$ to denote the scalar product on $M$ if there
    are no confusions.\\\indent  The second fundamental form in direction
$e_\alpha$ is denoted by \indent $$h_{\alpha ij}(x)= -\langle
e_{\alpha},\nabla_{i}\nabla_{j}F\rangle$$ and   the norm of the
second fundamental form by
$$|A|^{2}= g^{ij}g^{kl}h_{\alpha ik}h_{\alpha lj}.$$
The mean curvature on $M$ in direction $e_\alpha$ is given by
$$H_{\alpha}=g^{ij}h_{\alpha ij}.$$   Let $R_{ijkl}$ denote the curvature tensor
and  $R^{\bot}_{\beta\alpha jk}$ denote the normal curvature tensor,
and recall Ricci equation and Gauss equation for the submannifold of
Euclid space \beqn R^{\bot}_{ \alpha \beta ij}= h_{\alpha
ik}h_{\beta jk} - h_{\alpha jk}h_{\beta ik}\eeqn  and  \beqn
R_{ijkl} = h_{\alpha ik}h_{\alpha jl} - h_{\alpha jk}h_{\alpha
il}.\eeqn Of course, $R^{\bot}$ is zero for hypersurface. Also, we
can write Weigarten equation and Codazzi equation \beqn
\nabla_{i}e_{\beta} = h_{\beta i}^l \nabla_lF +
C^{\gamma}_{i\beta}e_{\gamma}\eeqn and \beqn h_{\alpha ik,j} =
h_{\alpha ij,k}\eeqn where $C^{\gamma}_{i\beta}$ is the connection
coefficient of normal connection and
$C^{\gamma}_{i\beta}=-C^{\beta}_{i\gamma}.$  Besides,
we will use the following basic facts.  \\
 {\bf Proposition
2.1}[14,18]. {\sl For any hypersurfece $M$ in $R^{n+1}$, we have
\beqn
\nabla_i\nabla_{j}F &=& - h_{ij}\nu, \\
\nabla_{i}\nu &= &h_{i}^{l}\nabla_{l}F, \\
\nabla_{k} h_{ij}& =& \nabla_{j}h_{ik}, \\
\nabla_i\nabla_{j}H &=& \triangle{h_{ij}}-H h_{i}^{l}h_{lj}+
|A|^{2}h_{ij},\\
2h^{ij}\nabla_i\nabla_{j}H &=&
\triangle{|A|^{2}}-2|\nabla{A}|^2-2Z,\eeqn where  $\nu$ is the outer
 normal vector of  $M$ ,
$C=g^{ij}g^{kl}g^{st}h_{ik}h_{sj}h_{lt} =tr(A^3)$, and
 $Z=HC-|A|^{4}.$}\\
 {\bf Proposition
2.2}[12][18]. {\sl  Suppose flow (1.1)holds true for $t\in [0,T)$
with $T\leq\infty$, then we have the following equations in [0,T):
\beqn \frac{dg_{ij}}{dt}&=& -2f_{\alpha}h_{\alpha ij},\\
\frac{dg^{ij}}{dt} &=& 2f_{\alpha}h_{\alpha kl}g^{ik}g^{jl},\eeqn
and letting $f\equiv f_{\alpha}=f_1$ for a hypersurface, we have
\beqn \frac{dh_{ij}}{dt} &=& \nabla_i\nabla_{j}f - fh_i^lh_{lj},\\
\frac{d|A|^{2}}{dt} &=& 2h^{ij}\nabla_i\nabla_{j}f + 2fC. \eeqn}
\indent  The following theorem for  short time existence
of (1.1) is well-known due to the theory of PDE and the technique of De Turk[19].\\
 \indent {\bf Theorem 2.3}[19] {\sl The flow (1.1) is a system of qusilinear parabolic
 equations and there exists a maximal time $0<T\leq \infty$ such
 that (1.1) admits a smooth solution on [0,T).}\\
\section{Hypersurfaces case}\setcounter{equation}{00}
In this section, we will prove theorems 1.2 and 1.4 and corollary
1.5. The key step is to  calculate the
evolution equations of $|A|^{2}.$  \\
 {\bf Proposition 3.1} {\sl
Suppose flow $(1.1)^{'}$ holds true for $t\in [0,T)$ with
$T\leq\infty$, then we have the following equations in [0,T):
\indent\beqn\frac{d|A|^{2}}{dt} &=&\triangle{|A|^{2}}-2|\nabla{A}|^2
+ 2|A|^{4} -
 2h^{ij}(\nabla_{i}\overline{\nabla}^{2}{\omega}){(\nabla_{j}{F},\nu)}
\nonumber\\&+&
2|A|^{2}{\overline{\nabla}^{2}}{\omega}(\nu,\nu)-4h^{ij}h_{j}^{l}\overline{\nabla}^{2}{\omega}
{(\nabla_{i}{F},\nabla_{l}{F})}-\langle\overline{\nabla}\omega,\nabla|A|^{2}\rangle.\eeqn}
{\bf Proof:} By (2.13), we have
$$\frac{d|A|^{2}}{dt} = 2h^{ij}\nabla_i\nabla_{j}f +
2fC.$$ By the notation $f\equiv f_{\alpha}=f_1$ and using (2.9), we
have \beqn 2h^{ij}\nabla_i\nabla_{j}f &=& 2h^{ij}\nabla_i\nabla_{j}H
-2h^{ij}\nabla_{i}\nabla_{j}\langle\overline{\nabla}\omega,\nu\rangle\nonumber
\\ &=& \triangle{|A|^{2}}-2|\nabla{A}|^2-2Z\ -
2h^{ij}\nabla_{i}\nabla_{j}\langle\overline{\nabla}\omega,
\nu\rangle .\eeqn   It follows from  (2.6) that
\beqs\nabla_{i}\nabla_{j}\langle\overline{\nabla}\omega,\nu\rangle
&=& \nabla_{i}(\langle\nabla_{j}\overline{\nabla}\omega,\nu\rangle +
\langle\overline{\nabla}\omega,\nabla_{j}\nu\rangle)\\
&=& \nabla_{i}(\langle\nabla_{j}\overline{\nabla}\omega,\nu\rangle +
h_{j}^{l}\langle\overline{\nabla}\omega, \nabla_{l}F\rangle )\\
&=& \langle\nabla_{i}\nabla_{j}\overline{\nabla}\omega,\nu\rangle +
h_{i}^{l}\langle\nabla_{j}\overline{\nabla}\omega,
\nabla_{l}F\rangle \\
&+&¡¡ h_{i}^{l}\langle\nabla_{j}\overline{\nabla}\omega,
\nabla_{l}F\rangle +
h_{j}^{l}\langle\overline{\nabla}\omega,\nabla_{i}\nabla_{l}F\rangle
+
\nabla_{i}h_{j}^{l}\langle\overline{\nabla}\omega,\nabla_{l}F\rangle.\eeqs
Using (2.5) and (2.7), we obtain   \indent
\beqn\nabla_{i}\nabla_{j}\langle\overline{\nabla}\omega,\nu\rangle
&=& \langle\nabla_{i}\nabla_{j}\overline{\nabla}\omega,\nu\rangle +
h_{i}^{l}\langle\nabla_{j}\overline{\nabla}\omega,
\nabla_{l}F\rangle \nonumber\\&+&
h_{i}^{l}\langle\nabla_{j}\overline{\nabla}\omega,
\nabla_{l}F\rangle + \langle\overline{\nabla}\omega,\nabla
h_{ij}\rangle -
h_{j}^{l}h_{il}\langle\overline{\nabla}\omega,\nu\rangle,\eeqn which
implies
\beqs2h^{ij}\nabla_{i}\nabla_{j}\langle\overline{\nabla}\omega,
\nu\rangle &=&
2h^{ij}\langle\nabla_{i}\nabla_{j}\overline{\nabla}\omega,\nu\rangle
+ 4h^{ij}h_{j}^{l}\langle\nabla_{i}\overline{\nabla}\omega,
\nabla_{l}F\rangle \\ &+&
\langle\overline{\nabla}\omega,\nabla|A|^{2}\rangle -
2h^{ij}h_{j}^{l}h_{il}\langle\overline{\nabla}\omega,\nu\rangle.\eeqs
This, together with (3.2),(2.13)and the definitions of $f$,C and Z, gives  \\
$$\frac{d|A|^{2}}{dt}
=\triangle{|A|^{2}}-2|\nabla{A}|^2 + 2|A|^{4} -
 2h^{ij}\langle\nabla_{i}\nabla_{j}\overline{\nabla}\omega,\nu\rangle
 $$
 \begin{equation}-4h^{ij}h_{j}^{l}\langle\nabla_{i}\overline{\nabla}\omega,
\nabla_{l}F\rangle -
\langle\overline{\nabla}\omega,\nabla|A|^{2}\rangle.\end{equation}
But \beqn
\langle\nabla_{i}\nabla_{j}\overline{\nabla}\omega,\nu\rangle &=&
\nabla_{i}(\langle\nabla_{j}{\overline{\nabla}{\omega}},\nu\rangle)
-\langle\nabla_{j}{\overline{\nabla}{\omega}},\nabla_{i}\nu\rangle \nonumber\\
&=& \nabla_{i}(\overline{\nabla}^{2}\omega(\nabla_{j}{F},\nu))-
\langle\nabla_{j}{\overline{\nabla}{\omega}},\nabla_{i}\nu\rangle\nonumber\\
&=&
(\nabla_{i}\overline{\nabla}^{2}{\omega}){(\nabla_{j}{F},\nu)}+\overline{\nabla}^{2}{\omega}(\nabla_{i}\nabla_{i}{F},
\nu)+
\langle\nabla_{j}{\overline{\nabla}{\omega}},\nabla_{i}\nu\rangle-
\langle\nabla_{j}{\overline{\nabla}{\omega}},\nabla_{i}\nu\rangle\nonumber\\
&=& (\nabla_{i}\overline{\nabla}^{2}{\omega}){(\nabla_{j}{F},\nu)}-
h_{ij}{\overline{\nabla}^{2}}{\omega}(\nu,\nu),\eeqn where we have
used (2.5) for the last equation.   Insert this
equality to (3.14) we can get the desired equality (3.1). \\
 {\bf Proposition 3.2} {\sl With the same assumption as in Proposition 3.1, the second
fundamental form satisfies the following evolution equation for
tensor in $[0, T):$ \beqn\frac{dh_{ij}}{dt} &=&\triangle{h_{ij}} -2H
h_{i}^{l}h_{lj}+ |A|^{2}h_{ij} -
 (\nabla_{i}\overline{\nabla}^{2}{\omega}){(\nabla_{j}{F},\nu)}\nonumber\\&+&
h_{ij}{\overline{\nabla}^{2}}{\omega}(\nu,\nu)-
h_{i}^{l}\overline{\nabla}^{2}{\omega}
{(\nabla_{j}{F},\nabla_{l}{F})} -
h_{i}^{l}\overline{\nabla}^{2}{\omega}
{(\nabla_{j}{F},\nabla_{l}{F})}\nonumber\\&-&
\langle\overline{\nabla}\omega,\nabla h_{ij}\rangle  +
2h_{j}^{l}h_{il}\langle\overline{\nabla}\omega,\nu\rangle.\eeqn}
{\bf
Proof:} It is a combination of (2.8),(2.12),(3.3) and (3.5).\\
{\bf Proof of Theorem 1.4:} Applying the maximum principle for
tensor [20] to equation (3.6), we see that
  the surface $M_t$ is
always convex along the flow if  $\overline{\nabla}^3{\omega}\equiv 0$ and
 $M_{0}$ is convex.\\
 \indent{\bf Lemma 3.3} {\sl Suppose that $M_{t}$ is the solution of
 $(1.1)^{'}$
 on [0,T) and  the assumptions (1), (2) and $(3)^{'}$ in theorem 1.2 are
 satisfied.
 Then $|A|^{2}<C$  for all $t\in  [0,T).$} \\
 {\bf Proof : } Taking a  local orthonomal basis $e_{i}$ on $M_t,
  i=1\ldots
n$,
 by (3.1) we have \beqs\frac{d|A|^{2}}{dt} & \leq &
\triangle{|A|^{2}}+ 2|A|^{4} -
 2h_{ij}({\nabla}_{i}\overline{\nabla}^{2}{\omega}){(e_{j},\nu)}
\nonumber\\&+&
2|A|^{2}{\overline{\nabla}^{2}}{\omega}(\nu,\nu)-4h_{ij}h_{jl}\overline{\nabla}^{2}{\omega}
{(e_{i},e_{l})}-\langle\overline{\nabla}\omega,\nabla|A|^{2}\rangle.\eeqs
 This, together with assumption (1), implies
\beqn\frac{d|A|^{2}}{dt} & \leq & \triangle{|A|^{2}}+ 2|A|^{4} +
 2|A|C_3 \nonumber\\&+& 2|A|^{2}\overline{\lambda}
-4h_{ij}h_{jl}\overline{\nabla}^{2}{\omega}
{(e_{i},e_{l})}-\langle\overline{\nabla}\omega,\nabla|A|^{2}\rangle.\eeqn
Next, we estimate $-4h_{ij}h_{jl}\overline{\nabla}^{2}{\omega}
{(e_{i},e_{l})}.$ Since
 \beqs
-4h_{ij}h_{jl}\overline{\nabla}^{2}{\omega} {(e_{i},e_{l})} &=
&-4h_{ij}h_{jl}(\overline{\lambda}E(e_{i},e_{l})
+\overline{\nabla}^{2}{\omega}
{(e_{i},e_{l})}-\overline{\lambda}E(e_{i},e_{l})) \nonumber\\&=&
-4|A|^{2}\overline{\lambda}-4h_{ij}h_{jl}(\overline{\nabla}^{2}\omega-\overline{\lambda}E)(e_{i},e_{l})\eeqs
 where $E$ is the unit matrix, we have  \beqs -4h_{ij}h_{jl}\overline{\nabla}^{2}{\omega}
{(e_{i},e_{l})} &\leq& -4|A|^{2}\overline{\lambda}
+4|h_{ij}h_{jl}||(\overline{\nabla}^{2}\omega-\overline{\lambda}E)(e_{i},e_{l})|
\nonumber \\ & \leq & -4|A|^{2}\overline{\lambda}+
4(\overline{\lambda}-\underline{\lambda})|A|^{2}\nonumber\\
&=& -4\underline{\lambda}|A|^{2}.\eeqs Therefore, (3.7) becomes
\beqn\frac{d|A|^{2}}{dt}  \leq  \triangle{|A|^{2}}+ 2|A|^{4} +
 2|A|C_3 + 2(\overline{\lambda}-2\underline{\lambda})|A|^{2}
-\langle\overline{\nabla}\omega,\nabla|A|^{2}\rangle.\eeqn Now by
the assumptions and Hamilton's maximum principle one can easily
obtain that $|A|^{2} <C$ for all  time.  Otherwise, we can choose
the first time $t_0$ such that $a(t_0)=C,$ where $a(t)\equiv
\max_{M_t} |A|^2.$ Then there exists a time $t_1<t_0$ such that
 $(\sqrt{C}-\delta )^2<a(t)<C$ for
$t\in [t_1, t_0)$ and $a(t_1)<C.$ Hence
$$2a^2(t)+2\sqrt{a(t)}C_3+2
(\overline{\lambda}-2\underline{\lambda})a(t)\leq 0, \forall t\in
[t_1, t_0)$$ by assumption $(3)'$. Therefore, applying  Hamilton's
maximum principle [20] to equation (3.8), we have $a(t)\leq
a(t_1)<C$ for all $t\in [t_1, t_0].$
This contradicts $a(t_0)=C.$ \\
{\bf Remark 3.4}  We would like to say that the convex condition on
$\omega$ (as in assumption (1)) and the small condition of the
initial $|A|^2$
 (as in assumptions (2) and $(3)'$) are necessary.  If  $\overline{\nabla}\omega=cx$ with  either $c<0$,
 or $c>0$ and $|A|^2>c$ on $M_0,$
  we have proved in
[17] that $|A|^{2}$ must blow up in finite time and the flow exists
only in
finite time.\\

 \indent {\bf Proof of Theorem 1.2:}   From  Lemma 3.3 we see that
 $|A|^{2}$  are bounded uniformly   if  assumptions (1)-(3) are satisfied.
 Thus, if we can prove that
$|\nabla^{m}A|^{2}\leq C_{m}$ is bounded when $t\longrightarrow T$,
then by a well-known theorem of partial differential equations the
flow (1.1) can be extended to [0,T+$\varepsilon$) for some small
$\varepsilon>0,$ where $T<\infty $ is the maximal time interval for
which $(1.1)'$ has a smooth solution.  This concludes
that the maximum time interval must be $[0, \infty ).$  \\
\indent To estimate $|\nabla A|^{2}$ , the boundedness of
$|\overline{\nabla}^{4}\omega|$ is necessary but is not enough,
because we  want to calculate the time derivative of
$\Gamma^{k}_{ij}$. As we know that connection is not a tensor, but
the difference of two connection is tensor, so is
$\frac{d\Gamma^{k}_{ij}}{dt}.$  Taking normal coordinate and using
(2.10), we have
\beqs\frac{d\Gamma^{k}_{ij}}{dt}&=&\frac{1}{2}\frac{d}{dt}(g^{lk}(\frac{\partial
g_{il}}{\partial x_j}+\frac{\partial g_{jl}}{\partial
x_{i}}-\frac{\partial g_{ij}}{\partial x_l})) \nonumber\\&=&
\frac{1}{2} (g^{lk}(\frac{\partial}{\partial
x_i}(\frac{d}{dt}g_{jl})+\frac{\partial}{\partial
x_j}(\frac{d}{dt}g_{il})-\frac{\partial}{\partial
x_l}(\frac{d}{dt}g_{ij}))
\nonumber\\&=&-g^{lk}(\frac{\partial}{\partial
x_i}(fh_{jl})+\frac{\partial}{\partial
x_j}(fh_{il})-\frac{\partial}{\partial x_l}(fh_{ij})).\eeqs  Noting
that
$\partial_{i}f=\partial_{i}H-\overline{\nabla}^{2}\omega(\partial_{i},\nu)-h_{il}\langle\overline{\nabla}\omega,
\partial_{l}\rangle$ and repeating the arguments of Huisken [14] we obtain  the following result.\\
{\bf Lemma 3.5}  {\sl Suppose that $M_{t}$ is the solution of
$(1.1)^{'}$
 on $ [0,T)$ for $T<\infty.$  If  assumptions (1),(2) and $(3)^{'}$ of theorem 1.2  are satisfied and
 $|\overline{\nabla}^{i}\omega|$ for $i=1,\ldots,m$ is uniformly bounded on $M_{t}$, then
$|\nabla^{m-3} A|^2$ is uniformly bounded on $M_{t}$.}\\
 Using lemma   3.5,   we have completed the proof of theorem 1.2.\\

{\bf Proof of corollary 1.5:} \indent   For the special case,
$\omega=\frac{1}{2}c_{1}x_{1}^{2}+ \ldots +
 \frac{1}{2}c_{n+1}x_{n+1}^2,$
 we have that\\\indent \indent  \indent  $\overline{\nabla}\omega=(c_{1}x_1,\ldots, c_{n+1}x_{n+1})$,
    $\overline{\nabla}^{2}\omega=(c_{i}\delta_{ij})$,   and
    $\overline{\nabla}^{3}\omega=0$ .\\
Let $M=\max{c_i}$, and  $m=\min{c_i}$ . Applying lemma 3.3 we get if
$M<2m$ and $|A|^{2}< 2m-M$ on $M_0$ , $|A|^{2}< 2m-M$ as long as
flow $(1.1)'$ exists. To get the long-time existence we have to get
the higher derivative estimate of $|A|^{2}.$  But lemma 3.5 can not
be applied directly, because $|\overline{\nabla}\omega |$ may turn
to be infinite if the surface expands to infinity. However, we can
prove that the surface will not expand to infinity in finite time as
follows. For this
purpose, we need  a theorem of [18].\\
{\bf Lemma 3.6}[18] {\sl Let F be a smooth immersed solution of
$(1.1)^{'}$ and $\widetilde{F}$ be an immersed solution of this
evolution equation. If $\widetilde{F}$ is contained in a connected
component of $R^{n+1}\setminus F$ or in the closure of such a
component at the beginning of the evolution, then this remains
during the evolution.}
\\
Since $|A|^{2}\leq 2m-M$ on $M_0$ , we will prove
that if the initial surface is sphere, the sphere will expand to
infinity as
$t\rightarrow\infty$ . \\
{\bf Lemma 3.7} {\sl Suppose that $M_0=S^{n}(R)$ is the initial
surface of the flow $(1.1)^{'}$ and $\omega ,$ $m$ and $M$ are as
above. Let  $s(t):=\frac{1}{2}|F_t|^2$ where$F_t$ is the position
vector of $M_t.$  If $|A|^{2}< 2m-M$ on $M_0$ , then $C_0\equiv
(2ms(0)-n)>0$ and $ s\geq \frac{n+C_0}{2m}e^{2mt}$ for all $t>0.$}

 {\bf
Proof :} Note that \beqs
\frac{ds}{dt}&=&\langle\frac{dF}{dt},F\rangle=-(H-\langle
\overline{\nabla}\omega,\nu\rangle)\langle F, \nu\rangle\nonumber\\
&=& -n+\langle\overline{\nabla}\omega,\nu\rangle\langle F,
\nu\rangle .\eeqs Since on the spheres $\nu=\frac{1}{|F|}F,$   we
have
\beqs\langle\overline{\nabla}\omega,\nu\rangle&=&\frac{1}{|F|}\langle\overline{\nabla}\omega,F\rangle\nonumber\\
&=& \frac{1}{|F|}(c_{1}F_{1}^2+\ldots +
c_{n+1}F_{n+1}^2)\nonumber\\
&\geq &\frac{1}{|F|}m|F|^2 = m|F|.\eeqs Hence, \beqn
\frac{ds}{dt}\geq -n + 2ms.\eeqn Therefore, $ s\geq
\frac{n+C_0}{2m}e^{2mt}$ for all $t>0$ if $C_0>0.$ Now by the
initial condition, we have \beqs 2m-M>|A|^{2} = \frac{1}{n}H^2 =
\frac{1}{n}\frac{n^2}{|F|^2} = \frac{n}{2s(0)},\eeqs which implies
 $2s>\frac{n}{2m-M}$ and $C_0>0.$  In this way, we have completed
 the proof of lemma 3.7.

\indent Finally, Lemma 3.6 and lemma 3.7 imply $M_t$ will not expand
to infinity in finite time. This, together with   the above
discussions, finishes the proof of corollary 1.5.

 \section{Higher co-dimension case}\setcounter{equation}{00}
 \indent In this section, we will prove theorem 1.1.  As the  hypersurface case,
 the key step is to derive the
 evolution equation of $|A|^{2}$.  For this purpose, we want to calculate
  the evolution equation of the
 second fundamental form tensor. In the following, for $x\in M^n$ we take orthonormal
 basis $e_1, \cdots , e_n, e_{n+1}, \cdots, e_{n+k}$ of $R^{n+k}$ such that $\{e_1, \cdots , e_n\}$
 is the basis of $T_xM^n$ and $\{e_{n+1}, \cdots, e_{n+k}\}$ (denoted by $\{e_{\alpha}\}$) is
  the unit normal vector.\\
{\bf Proposition 4.1} {\sl Suppose flow (1.1) holds true for $t\in
[0,T)$ with $T\leq\infty$, then we have the following equations in
[0,T):\beqn\frac{dh_{\alpha ij}}{dt} - \triangle h_{\alpha ij}&=&
-H_{\beta}h_{\beta il}h_{\alpha
jl}-(\nabla_{j}\overline{\nabla}^{2}{\omega}){(e_i,e_{\alpha})}+
h_{\beta
ij}{\overline{\nabla}^{2}}{\omega}(e_{\beta},e_{\alpha})\nonumber\\
&-& h_{\alpha kj}{\overline{\nabla}^{2}}{\omega}(e_{i},e_{k})-
h_{\alpha ik}{\overline{\nabla}^{2}}{\omega}(e_{j},e_{k}) + h_{\beta
ij}\langle e_{\beta}, \frac{de_\alpha}{dt}\rangle\nonumber\\
&+& \langle \overline{\nabla}\omega, e_\beta \rangle(h_{\beta
ik}h_{\alpha jk}+h_{\beta jk}h_{\alpha ik}) - \langle
\overline{\nabla}\omega, \nabla {h_{\alpha ij}}\rangle \nonumber\\
&-& h_{\alpha im}(h_{\gamma mj}h_{\gamma kk} - h_{\gamma
mk}h_{\gamma kj}) - h_{\alpha mk}(h_{\gamma mj}h_{\gamma
ik}-h_{\gamma mk}h_{\gamma ij}) \nonumber\\
&-& h_{\beta ik}(-h_{\beta km}h_{\alpha jm}+h_{\beta jm}h_{\alpha
km}).\eeqn }{\bf Proof: } For both sides are tensor, we calculate in
normal coordinate. Since $\nabla_{j}\nabla_{i}F= - h_{\alpha
ij}e_{\alpha}$ , then by (1.1) we have \beqn \frac{dh_{\alpha
ij}}{dt} &=& -\frac{d}{dt}\langle
\nabla_{j}\nabla_{i}F,e_{\alpha}\rangle \nonumber\\&=& - \langle
\nabla_{j}\nabla_{i}(-H_{\beta}e_{\beta}+\omega_{\beta}e_{\beta}),e_{\alpha}\rangle
- \langle\nabla_{j}\nabla_{i}F,\frac{d e_{\alpha}}{dt}\rangle
\nonumber\\&=& \langle
\nabla_{j}\nabla_{i}(H_{\beta}e_{\beta}),e_{\alpha}\rangle - \langle
\nabla_{j}\nabla_{i}(\omega_{\beta}e_{\beta}),e_{\alpha}\rangle +
 h_{\beta
ij}\langle e_{\beta}, \frac{de_\alpha}{dt}\rangle .\eeqn By
Weigarten equation(2.3), we have \beqn \nabla_{j}\nabla_{i}e_{\beta}
&=& ({\nabla_{j} h_{\beta il}}) e_{l} - h_{\beta il} h_{\gamma
jl}e_{\gamma}+ (\nabla_{j}{C^{\gamma}_{i\beta}})e_{\gamma}+
C^{\gamma}_{i\beta}\nabla_{j}e_{\gamma}\nonumber\\
&=& ({\nabla_{j} h_{\beta il}}) e_{l} - h_{\beta il} h_{\gamma
jl}e_{\gamma}+
(\nabla_{j}{C^{\gamma}_{i\beta}})e_{\gamma}+{C^{\gamma}_{i\beta}}h_{\gamma
jl}e_{l} + {C^{\gamma}_{i\beta}}{C^{\eta}_{j\gamma}}e_{\eta}\nonumber\\
&=& h_{\beta il,j} e_{l} - h_{\beta il} h_{\gamma jl}e_{\gamma}+
(\nabla_{j}{C^{\gamma}_{i\beta}})e_{\gamma} +
{C^{\gamma}_{i\beta}}{C^{\eta}_{j\gamma}}e_{\eta}.\eeqn This,
together with (2.3), implies
\beqn\nabla_{j}\nabla_{i}(H_{\beta}e_{\beta})
&=&(\nabla_{j}\nabla_{i}H_{\beta})e_{\beta} +
(\nabla_{j}H_{\beta})\nabla_{i}e_{\beta} +
(\nabla_{i}H_{\beta})\nabla_{j}e_{\beta} +
H_{\beta}\nabla_{j}\nabla_{i}e_{\beta}\nonumber\\
&=&
(\nabla_{j}\nabla_{i}H_{\beta})e_{\beta}+(\nabla_{j}H_{\beta})h_{\beta
il}e_{l} + (\nabla_{j}H_{\beta}){C^{\gamma}_{i\beta}}e_{\gamma}
\nonumber\\&+& (\nabla_{i}H_{\beta})h_{\beta jl}e_{l} +
(\nabla_{i}H_{\beta}{C^{\gamma}_{j\beta}})e_{\gamma} +
H_{\beta}\nabla_{j}\nabla_{i}e_{\beta}.\eeqn Hence, \beqn \langle
\nabla_{j}\nabla_{i}(H_{\beta}e_{\beta}),e_{\alpha}\rangle
&=&\nabla_{j}\nabla_{i}H_{\alpha} +
\nabla_{j}H_{\beta}{C^{\alpha}_{i\beta}} +
\nabla_{i}H_{\beta}{C^{\alpha}_{j\beta}} \nonumber\\&-
&H_{\beta}h_{\alpha jl}h_{\beta il} +
H_{\beta}\nabla_{j}{C^{\alpha}_{i\beta}} +
H_{\beta}{C^{\gamma}_{i\beta}}{C^{\eta}_{j\gamma}}.\eeqn Note that
\beqn \sum_{k}h_{\alpha kk, ij}&=&\nabla_{j}\nabla_{i}H_{\alpha} +
\nabla_{j}H_{\beta}{C^{\alpha}_{i\beta}} +
\nabla_{i}H_{\beta}{C^{\alpha}_{j\beta}}\nonumber\\
&+& H_{\beta}\nabla_{j}{C^{\alpha}_{i\beta}} +
H_{\beta}{C^{\gamma}_{i\beta}}{C^{\eta}_{j\gamma}}- 2h_{\alpha
kl}\frac{\partial{\Gamma^{l}_{ik}}}{\partial {x_j}} \eeqn and the
last term of (4.6) is zero because
$\Gamma^{l}_{ik}=-\Gamma^{k}_{il}$. Then we use  (4.6) to rewrite
(4.5) as\beqn \langle
\nabla_{j}\nabla_{i}(H_{\beta}e_{\beta}),e_{\alpha}\rangle =
\sum_{k}h_{\alpha kk, ij}-H_{\beta}h_{\alpha jl}h_{\beta il}.\eeqn
Simon's Identity gives\beqn \sum_{k}h_{\alpha kk, ij}=\triangle
h_{\alpha ij} -( h_{ \beta ik}R^{\bot}_{\beta\alpha jk} + h_{ \alpha
mk}R_{mijk} + h_{ \alpha im}R_{mkjk}).\eeqn Putting (4.8) in (4.7)
and using (2.1) and (2.2), we obtain that \beqn
 \langle
\nabla_{j}\nabla_{i}(H_{\beta}e_{\beta}),e_{\alpha}\rangle &=&
\triangle h_{\alpha ij} - H_{\beta}h_{\alpha jl}h_{\beta
il}\nonumber\\ &-& h_{\alpha im}(h_{\gamma mj}h_{\gamma kk} -
h_{\gamma mk}h_{\gamma kj}}) - h_{\alpha mk}(h_{\gamma mj}h_{\gamma
ik}-h_{\gamma mk}h_{\gamma ij) \nonumber\\
&-& h_{\beta ik}(-h_{\beta km}h_{\alpha jm} + h_{\beta jm}h_{\alpha
km}).\eeqn Next,  we use (2.3) to calculate the term
$\nabla_{j}\nabla_{i}(\omega_{\beta}e_{\beta})$ in (4.2). Since
\beqn\nabla_{j}\nabla_{i}(\omega_{\beta}e_{\beta})&=&
\nabla_{j}\nabla_{i}(\overline{\nabla}\omega -
\langle\overline{\nabla}\omega,e_{k}\rangle e_{k})\nonumber\\
&=& \nabla_{j}\nabla_{i}\overline{\nabla}\omega - \nabla_{j}
(\overline{\nabla}^{2}\omega(e_i,e_k)e_k \nonumber\\&-& h_{\beta
ik}\langle\overline{\nabla}\omega,e_{\beta}\rangle e_{k} - h_{\beta
ik}\langle\overline{\nabla}\omega,e_{k}\rangle
e_{\beta})\nonumber\\&=&\nabla_{j}\nabla_{i}\overline{\nabla}\omega
- \nabla_{j}(\overline{\nabla}^{2}\omega(e_i,e_k))e_k + h_{\beta
jk}\overline{\nabla}^{2}\omega(e_i,e_k)e_{\beta}\nonumber\\&+&
\nabla_{j}(h_{\beta
ik}\langle\overline{\nabla}\omega,e_{\beta}\rangle)e_{k} - h_{\beta
ik}h_{\gamma jk}\langle\overline{\nabla}\omega,e_{\beta}\rangle
e_{\gamma}\nonumber\\&+& \nabla_{j}(h_{\beta
ik})\langle\overline{\nabla}\omega,e_{k}\rangle e_{\beta} + h_{\beta
ik}\overline{\nabla}^{2}\omega(e_j,e_k)e_{\beta}\nonumber\\ &-&
h_{\beta ik}h_{\gamma
jk}\langle\overline{\nabla}\omega,e_{\gamma}\rangle e_\beta
+h_{\beta ik}\langle\overline{\nabla}\omega,e_{k}\rangle
\nabla_{j}e_\beta , \eeqn we have \beqn
\langle\nabla_{j}\nabla_{i}(\omega_{\beta}e_{\beta}),e_{\alpha}\rangle&=&
\langle\nabla_{j}\nabla_{i}\overline{\nabla}\omega,e_{\alpha}\rangle
+ h_{\alpha jk}\overline{\nabla}^{2}\omega(e_i,e_k) + h_{\alpha
ik}\overline{\nabla}^{2}\omega(e_j,e_k)\nonumber\\
&-&\langle\overline{\nabla}\omega,e_{\beta}\rangle(h_{\beta
ik}h_{\alpha jk}+h_{\alpha ik}h_{\beta jk})\nonumber\\&+&
(\nabla_{j}h_{\alpha ik}+C^{\alpha}_{j\beta}h_{\beta
ik})\langle\overline{\nabla}\omega,e_{k}\rangle\nonumber\\&=&
\langle\nabla_{j}\nabla_{i}\overline{\nabla}\omega,e_{\alpha}\rangle
+ h_{\alpha jk}\overline{\nabla}^{2}\omega(e_i,e_k) + h_{\alpha
ik}\overline{\nabla}^{2}\omega(e_j,e_k)\nonumber\\
&-&\langle\overline{\nabla}\omega,e_{\beta}\rangle(h_{\beta
ik}h_{\alpha jk}+h_{\alpha ik}h_{\beta jk})+ h_{\alpha
ik,j}\langle\overline{\nabla}\omega,e_{k}\rangle .\eeqn Due to
Codazzi equation (2.4), we have\beqn
\langle\nabla_{j}\nabla_{i}(\omega_{\beta}e_{\beta}),e_{\alpha}\rangle
&=&
\langle\nabla_{j}\nabla_{i}\overline{\nabla}\omega,e_{\alpha}\rangle
+ h_{\alpha jk}\overline{\nabla}^{2}\omega(e_i,e_k) + h_{\alpha
ik}\overline{\nabla}^{2}\omega(e_j,e_k)\nonumber\\
&-&\langle\overline{\nabla}\omega,e_{\beta}\rangle(h_{\beta
ik}h_{\alpha jk}+h_{\alpha ik}h_{\beta jk})+
\langle\overline{\nabla}\omega,\nabla h_{\alpha
ij}\rangle\nonumber\\
&=&(\nabla_{j}\overline{\nabla}\omega^{2})(e_i,e_{\alpha})-h_{\beta
ij}\overline{\nabla}^2\omega (e_\alpha,e_{\beta})\nonumber\\&+
&h_{\alpha jk}\overline{\nabla}^{2}\omega(e_i,e_k) +h_{\alpha
ik}\overline{\nabla}^{2}\omega(e_j,e_k)\nonumber\\
&-&\langle\overline{\nabla}\omega,e_{\beta}\rangle(h_{\beta
ik}h_{\alpha jk}+h_{\alpha ik}h_{\beta jk})+
\langle\overline{\nabla}\omega,\nabla h_{\alpha ij}\rangle .\eeqn So
(4.1) follows from (4.2),(4.9) and (4.12).\\
{\bf Proposition 4.2} {\sl Suppose flow (1.1)holds true for $t\in
[0,T)$ with $T\leq\infty$, then we have the following equation in
[0,T): \beqn \frac{d|A|^{2}}{dt}&=&\triangle|A|^{2}-2|\nabla
A|^{2}-2h_{\alpha
ij}(\nabla_{j}\overline{\nabla}^2\omega)(e_i,e_{\alpha})\nonumber\\
&+&2h_{\alpha ij}h_{\beta
ij}\overline{\nabla}^{2}\omega(e_\alpha,e_\beta)-4h_{\alpha
ik}h_{\alpha
ij}\overline{\nabla}^{2}\omega(e_j,e_k)-\langle\overline{\nabla}\omega,\nabla
|A|^2\rangle\nonumber\\&+&2\sum_{\alpha,\gamma,i,m}(\sum_{k}h_{\alpha
ik}h_{\gamma mk}-h_{\alpha mk}h_{\gamma
ik})^{2}+2\sum_{i,j,k,m}(\sum_{\alpha}h_{\alpha ij}h_{\alpha
mk})^{2}\eeqn } {\bf Proof: } We calculate it in normal coordinate.
Because $|A|^{2}= g^{ij}g^{kl}h_{\alpha ik}h_{\alpha lj}$, then
\beqn\frac{d|A|^{2}}{dt}&=&2\frac{dg^{ik}}{dt}h_{\alpha ij}h_{\alpha
kj}+2\frac{dh_{\alpha ij}}{dt}h_{\alpha ij}.\eeqn Hence by (2.11)
(4.1) and (4.14), we have\beqs \frac{d|A|^{2}}{dt}&=&2h_{\alpha
ij}\triangle h_{\alpha ij}+4(H_{\beta}-\omega_{\beta})h_{\beta
ik}h_{\alpha ij}h_{\alpha kj}-2H_{\beta}h_{\alpha ij}h_{\beta
il}h_{\alpha
jl}\nonumber\\
&-&2h_{\alpha
ij}(\nabla_{j}\overline{\nabla}^{2}{\omega}){(e_i,e_{\alpha})}+2h_{\alpha
ij} h_{\beta
ij}{\overline{\nabla}^{2}}{\omega}(e_{\beta},e_{\alpha})\nonumber\\
&-&4h_{\alpha ij}h_{\alpha
kj}{\overline{\nabla}^{2}}{\omega}(e_{i},e_{k}) +  2h_{\alpha
ij}h_{\beta ij}\langle e_{\beta}, \frac{de_\alpha}{dt}\rangle
+4h_{\alpha ij} \omega_\beta h_{\beta ik}h_{\alpha jk}\nonumber\\&-
& \langle \overline{\nabla}\omega, \nabla |A|^{2}\rangle -
2h_{\alpha ij}h_{\alpha im}h_{\gamma mj}H_{\gamma } + 2h_{\alpha
ij}h_{\alpha im}h_{\gamma mk}h_{\gamma kj} \nonumber\\ &-&
2h_{\alpha ij}h_{\alpha mk}(h_{\gamma mj}h_{\gamma ik}-h_{\gamma
mk}h_{\gamma ij}) -2h_{\alpha ij}h_{\beta ik}(h_{\beta lj}h_{\alpha
lk}- h_{\beta lk}h_{\alpha lj}).\eeqs Observing that $2h_{\alpha
ij}h_{\beta ij}\langle e_{\beta}, \frac{de_\alpha}{dt}\rangle$ is
zero by symmetry and $2h_{\alpha ij}\triangle h_{\alpha
ij}=\triangle |A|^{2}-2|\nabla A|^2$ , we have \beqs
\frac{d|A|^{2}}{dt}&=&\triangle |A|^{2}-2|\nabla A|^2 -2h_{\alpha
ij}(\overline{\nabla}_{j}\overline{\nabla}^{2}{\omega}){(e_i,e_{\alpha})}+2h_{\alpha
ij} h_{\beta
ij}{\overline{\nabla}^{2}}{\omega}(e_{\beta},e_{\alpha})\nonumber\\
&-&4h_{\alpha ij}h_{\alpha
kj}{\overline{\nabla}^{2}}{\omega}(e_{i},e_{k}) -  \langle
\overline{\nabla}\omega, \nabla |A|^{2}\rangle  + 2h_{\alpha
ij}h_{\alpha im}h_{\gamma mk}h_{\gamma kj} \nonumber\\ &-&
2h_{\alpha ij}h_{\alpha mk}(h_{\gamma mj}h_{\gamma ik}-h_{\gamma
mk}h_{\gamma ij}) -2h_{\alpha ij}h_{\beta ik}(h_{\beta lj}h_{\alpha
lk}- h_{\beta lk}h_{\alpha lj}).\eeqs But the last three terms can
be calculate as follows: \beqn&\indent&2h_{\alpha ij}h_{\alpha
im}h_{\gamma mk}h_{\gamma kj} - 2h_{\alpha ij}h_{\alpha mk}h_{\gamma
mj}h_{\gamma ik}\nonumber\\&+&2h_{\alpha ij}h_{\alpha mk}h_{\gamma
mk}h_{\gamma ij} -2h_{\alpha ij}h_{\beta ik}(h_{\beta lj}h_{\alpha
lk}- h_{\beta lk}h_{\alpha lj})\nonumber\\&=&4h_{\alpha ij}h_{\alpha
im}h_{\gamma mk}h_{\gamma kj}-4h_{\alpha ij}h_{\alpha mk}h_{\gamma
mj}h_{\gamma ik}+2h_{\alpha ij}h_{\gamma mk}h_{\alpha mk}h_{\gamma
ij}.\eeqn Since \beqn &\indent& 2h_{\alpha ij}h_{\alpha im}h_{\gamma
mk}h_{\gamma kj}-2h_{\alpha ij}h_{\alpha mk}h_{\gamma mj}h_{\gamma
ik}\nonumber\\&=& 2h_{\alpha ij}h_{\alpha ik}h_{\gamma mk}h_{\gamma
mj}-2h_{\alpha ij}h_{\alpha mk}h_{\gamma mj}h_{\gamma
ik}\nonumber\\&=& 2h_{\alpha ij}h_{\gamma mj}(h_{\alpha ik}h_{\gamma
mk}-h_{\alpha mk}h_{\gamma ik})\nonumber\\&=& h_{\alpha ij}h_{\gamma
mj}(h_{\alpha ik}h_{\gamma mk}-h_{\alpha mk}h_{\gamma ik})+
h_{\alpha mj}h_{\gamma ij}(h_{\alpha mk}h_{\gamma
ik}-h_{\alpha mk}h_{\alpha ik})\nonumber\\
&=& \sum_{\alpha,\gamma,i,m}(\sum_{k}h_{\alpha ik}h_{\gamma
mk}-h_{\alpha mk}h_{\gamma ik})^{2}\eeqn and \beqn2h_{\alpha
ij}h_{\gamma mk}h_{\alpha mk}h_{\gamma
ij}=2\sum_{i,j,k,m}(\sum_{\alpha}h_{\alpha ij}h_{\alpha
mk})^{2},\eeqn  we have proved the proposition.\\
{\bf Lemma 4.3}  {\sl Suppose that $M_{t}$ is the solution of (1.1)
 on [0,T), and the assumptions (1), (2) and (3) of theorem 1.1 are
 satisfied,
 then $|A|^{2}\leq C $  on $M_t$ for all $ t\in [0,T).$}

{\bf Proof: } The proof  is almost the same as  that of lemma 3.3 in
the case
 of hypersurfaces.  It follows from  Schwartz inequality
 that \beqs2\sum_{i,j,k,m}(\sum_{\alpha}h_{\alpha ij}h_{\alpha
mk})^{2}\leq 2\sum_{i,j,k,m}(\sum_{\alpha}h_{\alpha
ij}^{2})(\sum_{\alpha}h_{\alpha mk}^{2})=2|A|^4\eeqs and \beqs
\sum_{\alpha,\gamma,i,m}(\sum_{k}h_{\alpha ik}h_{\gamma
mk}-h_{\alpha mk}h_{\gamma ik})^{2}\leq
4\sum_{\alpha,\gamma,i,m}(\sum_{k}h_{\alpha ik}h_{\gamma
mk})^{2}\leq4|A|^4.\eeqs Consequently,  using the same technique
from (3.7) to (3.8)  we  obtain  \beqn
\frac{d|A|^{2}}{dt}&\leq&\triangle |A|^{2}- \langle
\overline{\nabla}\omega, \nabla |A|^{2}\rangle+10|A|^{4} +
 2|A|C_3 + 2(\overline{\lambda}-2\underline{\lambda})|A|^{2}.\eeqn
then the result follows by copying the arguments below (3.8).\\

{\bf Proof of theorem 1.1:}  Using Lemma 4.3 and repeating the proof
of theorem 1.2, one can easily prove theorem 1.1.

\end{document}